# Quantized control via locational optimization

Francesco Bullo* and Daniel Liberzon†

Coordinated Science Laboratory
University of Illinois at Urbana-Champaign
1308 W. Main St, Urbana, IL 61801, USA

November 12, 2018

## Abstract

This paper studies state quantization schemes for feedback stabilization of control systems with limited information. The focus is on designing the least destabilizing quantizer subject to a given information constraint. We explore several ways of measuring the destabilizing effect of a quantizer on the closed-loop system, including (but not limited to) the worst-case quantization error. In each case, we show how quantizer design can be naturally reduced to a version of the so-called multicenter problem from locational optimization. Algorithms for solving such problems are discussed. In particular, an iterative solver is developed for a novel weighted multicenter problem which most accurately represents the least destabilizing quantizer design.

## 1 Introduction

In this paper we study control systems whose state variables are quantized. We think of a quantizer as a device that converts a real-valued signal into a piecewise constant one taking a finite set of values. The recent papers [2, 7, 10] discuss various situations where this type of quantization arises and provide references to the literature. Mathematically, a quantizer can be described by a piecewise constant function $q : \mathcal{D} \subset \mathbb{R}^n \to \mathcal{Q}$, where $\mathcal{Q}$ is a finite subset of $\mathbb{R}^n$ with a fixed number of elements $N$. Here $n$ is the state dimension of a given system. We denote the elements of $\mathcal{Q}$ by $q_1, \ldots, q_N$ and refer to them as *quantization points*. The sets $W_i := \{x \in \mathcal{D} : q(x) = q_i\}$, $i \in \{1, \ldots, N\}$ associated with fixed values of the quantizer form a partition[1] of the domain $\mathcal{D}$ and are called *quantization regions*. Quantized values of the state represent a limited information flow from the system to a feedback controller: the state is not completely known to the controller, but it is only known which one of a fixed number of quantization regions contains the current state at each instant of time.

In the literature it is usually assumed that quantization regions are fixed in advance and have specific shapes, most often rectilinear. Here we are interested in the situation where the number of values $N$ of the quantizer is a given information constraint in the control problem, but the control designer has flexibility in choosing a specific configuration of quantization regions (whose shapes can in principle be arbitrary) and quantization points. While there has been some research on systems with quantization regions of arbitrary shapes [14, 12] and on the relationship between the choice of quantization regions and the behavior of the closed-loop system [7, 10], the general problem of determining the "best" quantizer for a particular control task such as feedback stabilization remains largely open.

---

*Email: bullo@uiuc.edu. Supported by DARPA/AFOSR Award F49620-02-1-0325.
†Email: liberzon@uiuc.edu. Supported by NSF ECS-0134115 CAR, NSF ECS-0114725, and DARPA/AFOSR MURI F49620-02-1-0325 grants.

[1]A collection $\{W_1, \ldots, W_N\}$ of subsets of $\mathcal{D}$ is a *partition* of $\mathcal{D}$ if the intersection between the relative interior of any two $W_i$ is empty and the union of all $W_i$ equals $\mathcal{D}$.

A feedback law which globally asymptotically stabilizes a given system in the absence of quantization will in general fail to provide global asymptotic stability of the closed-loop system that arises in the presence of state quantization. There are two phenomena which account for changes in the system's behavior caused by quantization. The first one is saturation: if the quantized signal is outside the range of the quantizer, then the quantization error is large, and the control law designed for the ideal case of no quantization leads to instability. The second one is deterioration of performance near the equilibrium: as the difference between the current and the desired values of the state becomes small, higher precision is required, and so in the presence of quantization errors asymptotic convergence is impossible. These phenomena manifest themselves in the existence of two nested invariant regions such that all trajectories of the quantized system starting in the bigger region approach the smaller one, while no further convergence guarantees can be given.

In Section 2 we explain how the destabilizing effect of a given quantizer can be measured. This *destabilization measure*, in conjunction with an arbitrary stabilizing feedback law and a corresponding Lyapunov function, can be used to determine an ultimate bound on solutions. One example of such a destabilization measure is the *worst-case quantization error* $\max_{x \in \mathcal{D}} |q(x) - x|$. However, it turns out that there exist other destabilization measures which are actually more suitable in the present context. Although the parameters of the control system are used in the stability analysis, the destabilization measure itself is a function of the quantization regions and quantization points only. The quantizer design problem then naturally reduces to an optimization problem which consists in minimizing such a measure over all quantizers satisfying the information constraint.

After casting quantizer design as an optimization problem, we proceed to explain how techniques from (optimal) facility location yield new insights into this problem as well as efficient algorithms for solving it. Facility location problems concern the location of a fixed number of facilities that provide service demanded by users; the objective is to minimize the average or maximal distance from sets of demand points to facilities. We focus here on settings continuous in the location of both the facilities and the demand points (both facilities and demand points take values in a continuum of points, e.g., in a polytope or ellipsoid). Facility location problems are surveyed in [5]. Computational geometric aspects in continuous facility location are discussed in [19, 18] and indirectly in textbooks on computational geometry [4]. Relevant background on computational geometric methods in locational optimization is provided in Section 3.

For example, a classical problem of interest in locational optimization is the so-called *multimedian problem*. It consists in choosing a collection of $N$ points $q_1, q_2, \ldots, q_N$ in a bounded region $\mathcal{D} \subset \mathbb{R}^n$ so as to minimize the quantity $E(\min_{i \in \{1,\ldots,N\}} |q_i - x|^2)$, where the expected value is computed with respect to some probability density function on $\mathcal{D}$ and $|\cdot|$ denotes the Euclidean norm. A solution of this problem is given by *centroidal Voronoi tessellations*; see [23, 6]. Within the context of quantization and information theory, the multimedian problem is known as the classic fixed-rate minimum-distortion quantizer design [9]. One of the early references on this problem is the work by Lloyd [13], who obtains optimality conditions and introduces a famous insightful algorithm. The multimedian problem is related to the problem of state moment stabilization of linear systems with limited data rate [16].

Since we are working in the deterministic setting, we will find that the problem relevant for our purposes is the *multicenter problem*, discussed in [23, 22]. This is a somewhat less frequently encountered variant of the multimedian problem, which is obtained by replacing the expected value by the worst-case value; it can also be stated as the problem of covering a given region with overlapping balls of minimal radius. The connection between the quantized control problem and the multicenter problem, although very natural, apparently has not been pursued before. In Section 3 we present a general formulation of the multicenter problem with weighting factors. We then discuss solutions of several versions of this problem which arise in the present context, all in terms of suitable Voronoi tessellations. We show how existing algorithms can handle some of the more apparent approaches, and then develop a novel algorithm which gives less conservative results. We note that while weighted multimedian problems are commonly encountered, our formulation is novel in that it introduces weighting factors in the multicenter problem.

Simulation studies that we are currently conducting (as well as existing studies of the related multimedian problem, such as [8]) indicate that by solving the quantized feedback stabilization problem with the help of locational optimization techniques, one may obtain quite interesting quantization patterns. In particular, for planar systems a typical quantization region is a hexagon. Consequently, hexagonal quantization regions are capable of achieving better performance than more traditional rectangular ones.

## 2 Quantization and stability

We assume that the stabilization problem in the absence of quantization has been solved, in the sense that a state feedback control law is known such that the origin is a globally asymptotically stable equilibrium point of the ideal closed-loop system. In the presence of quantization, we adopt the "certainty equivalence" control paradigm; namely, we let the same control law act on the quantized state $q(x)$. The problem under consideration is to characterize the destabilizing effect of a quantizer $q$, with the goal of obtaining an ultimate bound on solutions of the closed-loop system starting in a given bounded region. We first discuss this problem for general nonlinear systems and then develop more specific results for linear systems, moving from simpler but conservative to more complicated but sharper formulations.

### 2.1 Nonlinear systems

We start with the general situation where the process to be controlled is modeled by the nonlinear system
$$\dot{x} = f(x, u), \qquad x \in \mathbb{R}^n, \ u \in \mathbb{R}^m. \tag{1}$$
All vector fields and control laws are understood to be sufficiently regular (e.g., $\mathcal{C}^1$) so that existence and uniqueness of solutions are ensured. Suppose that some nominal static feedback law $u = k(x)$ is given (with minor changes, dynamic feedback laws can also be used). In the presence of state quantization, we consider the feedback law $u = k(q(x))$ and the corresponding closed-loop system
$$\dot{x} = f(x, k(q(x))) = f(x, k(x + e)) \tag{2}$$
where
$$e := q(x) - x$$
represents the *quantization error*.

Besides stabilizing the nominal system (1), the feedback law $k$ clearly must possess some robustness property with respect to the measurement error $e$. To this end we impose the following assumption: there exists a $\mathcal{C}^1$ function $V : \mathbb{R}^n \to \mathbb{R}$ such that for some class $\mathcal{K}_\infty$ functions[2] $\alpha_1, \alpha_2, \alpha_3, \rho$ and for all

---

[2]Recall that a function $\alpha : [0, \infty) \to [0, \infty)$ is said to be of *class* $\mathcal{K}$ if it is continuous, strictly increasing, and $\alpha(0) = 0$. If $\alpha \in \mathcal{K}$ is unbounded, then it is said to be of *class* $\mathcal{K}_\infty$. A function $\beta : [0, \infty) \times [0, \infty) \to [0, \infty)$ is said to be of *class* $\mathcal{KL}$ if $\beta(\cdot, t)$ is of class $\mathcal{K}$ for each fixed $t \geq 0$ and $\beta(s, t)$ decreases to 0 as $t \to \infty$ for each fixed $s \geq 0$.

$x, e \in \mathbb{R}^n$ we have
$$\alpha_1(|x|) \leq V(x) \leq \alpha_2(|x|)$$
and
$$|x| \geq \rho(|e|) \Rightarrow \frac{\partial V}{\partial x} f(x, k(x+e)) \leq -\alpha_3(|x|).$$

(Here and below, $|\cdot|$ denotes the standard Euclidean norm.) This amounts to the property that the control law $u = k(x)$ *input-to-state stabilizes* the closed-loop system with respect to the measurement error $e$ [20, 21]. The above assumption is rather restrictive and can be relaxed at the expense of obtaining weaker results (for linear systems, however, it is an automatic consequence of closed-loop asymptotic stability for $e \equiv 0$). There is also considerable research on designing control laws satisfying this assumption. These issues are discussed elsewhere [11, 12].

Pick a positive number $M$ and consider the ball $\mathcal{B}_M := \{x \in \mathbb{R}^n : |x| \leq M\}$. Consider the worst-case quantization error
$$\Delta := \max_{x \in \mathcal{B}_M} |e| \tag{3}$$
(this quantity is sometimes also referred to as the *sensitivity* of the quantizer). The following result is then straightforward to prove (see [12]).

**Lemma 1** *Assume that*
$$\alpha_1(M) > \alpha_2 \circ \rho(\Delta). \tag{4}$$
*Then the sets*
$$\mathcal{R}_1 := \{x \in \mathbb{R}^n : V(x) \leq \alpha_1(M)\} \tag{5}$$
*and*
$$\mathcal{R}_2 := \{x \in \mathbb{R}^n : V(x) \leq \alpha_2 \circ \rho(\Delta)\}$$
*are invariant regions for the system* (2). *Moreover, all solutions of* (2) *that start in the set $\mathcal{R}_1$ enter the smaller set $\mathcal{R}_2$ in finite time. An upper bound on this time is*
$$T = \frac{\alpha_1(M) - \alpha_2 \circ \rho(\Delta)}{\alpha_3 \circ \rho(\Delta)}. \tag{6}$$

This lemma implies, in particular, that all solutions starting in $\mathcal{R}_1$ at time $t = t_0$ satisfy the ultimate bound
$$|x(t)| \leq \alpha_1^{-1} \circ \alpha_2 \circ \rho(\Delta) \qquad \forall t \geq t_0 + T \tag{7}$$
with $T$ given by the formula (6). We regard the quantity $\Delta$ defined by (3) as a *destabilization measure* of the quantizer $q$. For given feedback law $k$ and Lyapunov function $V$, an ultimate bound on solutions can be described by a class $\mathcal{K}_\infty$ function of this measure as shown above. It is not hard to see that if the number $N$ of quantization regions is sufficiently large, then $\Delta$ can be made small enough for the inequality (4) to hold. Minimizing $\Delta$—and consequently the size of the attracting invariant region $\mathcal{R}_1$—over all possible choices of the quantizer $q$ corresponds to the following optimization problem:
$$\min_{\mathcal{Q}, \mathcal{W}} \max_{i \in \{1,\ldots,N\}} \max_{x \in W_i} |q_i - x| \tag{8}$$
where $\mathcal{Q} = \{q_1, \ldots, q_N\}$ is a set of quantization points and $\mathcal{W} = \{W_1, \ldots, W_N\}$ is a partition of $\mathcal{B}_M$ into quantization regions. (We could work with partitions of $\mathcal{R}_1$ rather than $\mathcal{B}_M$, but this requires the knowledge of $P$ and also may be less computationally feasible for non-quadratic $V$.) The optimization problem (8) is known as the *multicenter problem* in computational geometry; we defer its detailed discussion until Section 3.1.

## 2.2 Linear systems

We now specialize to the case when the process is described by the linear system

$$\dot{x} = Ax + Bu, \qquad x \in \mathbb{R}^n, \ u \in \mathbb{R}^m. \tag{9}$$

The linear system structure can be utilized to define a less conservative destabilization measure. Suppose that the system (9) is stabilizable, so that for some matrix $K$ the eigenvalues of $A + BK$ have negative real parts. Then there exists a unique positive definite symmetric matrix $P$ such that

$$(A + BK)^T P + P(A + BK) = -I. \tag{10}$$

We let $\lambda_{\min}(P)$ and $\lambda_{\max}(P)$ denote the smallest and the largest eigenvalue of $P$, respectively. We assume that $B \neq 0$ and $K \neq 0$; this is no loss of generality because the case of interest is when $A$ is not a stable matrix.

The quantized state feedback control law

$$u = Kq(x)$$

yields the closed-loop system

$$\dot{x} = Ax + BKq(x) = (A + BK)x + BKe \tag{11}$$

where $e := q(x) - x$ is the quantization error as before. The derivative of the function

$$V(x) := x^T Px$$

along solutions of the system (11) satisfies

$$\dot{V} = -x^T x + 2x^T PBKe \leq -|x|^2 + 2|x||PBKe|. \tag{12}$$

For an arbitrary small $\varepsilon > 0$, we can rewrite this as

$$\dot{V} = -\frac{\varepsilon}{1+\varepsilon}|x|^2 - \frac{1}{2(1+\varepsilon)}|x|^2 + 2(1+\varepsilon)|PBKe|^2 - \left(\frac{1}{\sqrt{2(1+\varepsilon)}}|x| - \sqrt{2(1+\varepsilon)}|PBKe|\right)^2$$

$$\leq -\frac{\varepsilon}{1+\varepsilon}|x|^2 - \frac{1}{2(1+\varepsilon)}\left(|x|^2 - \left(2(1+\varepsilon)|PBKe|\right)^2\right)$$

Therefore, we have

$$|x| \geq 2(1+\varepsilon)|PBKe| \quad \Rightarrow \quad \dot{V} \leq -\frac{\varepsilon}{1+\varepsilon}|x|^2. \tag{13}$$

Take a positive number $M$. There are two ways to proceed from here. One is to use the inequalities (13) and $|PBKe| \leq \|PBK\||e|$, where $\|\cdot\|$ stands for the induced matrix norm, to obtain

$$|x| \geq 2(1+\varepsilon)\|PBK\||e| \quad \Rightarrow \quad \dot{V} \leq -\frac{\varepsilon}{1+\varepsilon}|x|^2. \tag{14}$$

Consider the ellipsoids

$$\mathcal{R}_1 := \{x \in \mathbb{R}^n : x^T Px \leq \lambda_{\min}(P)M^2\} \tag{15}$$

and

$$\mathcal{R}_2 := \{x \in \mathbb{R}^n : x^T Px \leq \lambda_{\max}(P)4(1+\varepsilon)^2\|PBK\|^2\Delta^2\} \tag{16}$$

where $\Delta$ is the worst-case quantization error defined by (3). Then we have the following linear counterpart of Lemma 1.

**Lemma 2** *Assume that*

$$\lambda_{\min}(P)M^2 > \lambda_{\max}(P)4(1+\varepsilon)^2\|PBK\|^2\Delta^2. \tag{17}$$

*Then the ellipsoids $\mathcal{R}_1$ and $\mathcal{R}_2$ are invariant regions for the system (11). Moreover, all solutions of (11) that start in the ellipsoid $\mathcal{R}_1$ enter the smaller ellipsoid $\mathcal{R}_2$ in finite time. An upper bound on this time is*

$$T = \frac{\lambda_{\min}(P)M^2 - \lambda_{\max}(P)4(1+\varepsilon)^2\|PBK\|^2\Delta^2}{4\|PBK\|^2\Delta^2(1+\varepsilon)\varepsilon}. \tag{18}$$

Consequently, an ultimate bound on solutions starting in $\mathcal{R}_1$ at time $t = t_0$ is

$$|x(t)| \leq \sqrt{\frac{\lambda_{\max}(P)}{\lambda_{\min}(P)}} 2(1+\varepsilon)\|PBK\|\Delta \qquad \forall t \geq t_0 + T$$

with $T$ given by the formula (18). Decreasing $\varepsilon$ to 0, we see that solutions (asymptotically) approach the ellipsoid

$$\{x \in \mathbb{R}^n : x^T P x \leq \lambda_{\max}(P)4\|PBK\|^2\Delta^2\}.$$

Thus we still consider $\Delta$ as a destabilization measure. As in the nonlinear setting, this leads to the optimization problem (8). If $N$ is large enough, then $\Delta$ can be made small enough so that the inequality (17) holds (for a given feedback gain $K$).

Another approach is to work with (13) directly, avoiding the use of the induced norm $\|PBK\|$. Define

$$\Delta_{PBK} := \max_{x \in \mathcal{R}_1} |PBKe|.$$

The result of Lemma 2 still holds if $\|PBK\|^2\Delta^2$ is replaced by $\Delta_{PBK}^2$ everywhere in the statement of that lemma. This yields a less conservative ultimate bound and motivates the following optimization problem:

$$\min_{\mathcal{Q},\mathcal{W}} \max_{i \in \{1,\ldots,N\}} \max_{x \in W_i} |PBK(q_i - x)| \tag{19}$$

where $\mathcal{Q}$ is a set of quantization points as before and $\mathcal{W}$ is a partition of $\mathcal{R}_1$ into $N$ regions. This problem is in general lower-dimensional compared to (8) because the subspace $\ker(PBK)$ can be ignored. (Note that $PBK$ is a singular matrix whenever $m < n$.) Therefore, for the same $N$ the optimal value for this problem will be significantly lower that for (8). However, $\Delta_{PBK}$ is not really a destabilization measure in the sense used in this paper, because it depends on the feedback gain matrix $K$. While it gives better results for a fixed feedback law, quantizer design based on this destabilization measure needs to be redone if the feedback law is changed, and is not suitable for switching between several feedback laws. For these reasons we prefer to work with destabilization measures that are independent of a particular feedback law used.

**Remark 1** It is clear that in Lemma 2, the system's behavior is important only for $x \in \mathcal{R}_1 \setminus \mathcal{R}_2$. We can thus work with partitions of $\mathcal{R}_1 \setminus \mathcal{R}_2$—or of a spherical annulus containing this set—rather than with partitions of $\mathcal{R}_1$. This preserves or decreases the ultimate bound on solutions. □

**Remark 2** Lemma 2 suggests that among stabilizing state feedback gains $K$, the ones that provide smaller ultimate bounds for the solutions of the quantized system are those with smaller values of the induced matrix norm $\|PBK\|$, where $P$ is given by (10). In this regard, it might be interesting to observe the following: if the open-loop system $\dot{x} = Ax$ is not asymptotically stable, then for every stabilizing feedback gain $K$ and the corresponding positive definite symmetric matrix $P$ satisfying (10) we have
$$\|PBK\| \geq 1/2 \tag{20}$$
and the inequality is strict if $\dot{x} = Ax$ is unstable. To see this, use (12) and the definition of $e$ to write
$$\dot{V} \leq -|x|^2 \left(1 - 2\|PBK\|\frac{|q(x) - x|}{|x|}\right). \tag{21}$$
This formula will be used again several times in the sequel. Now, note that if $q$ is chosen to take the value 0 in a neighborhood of the origin, then the right-hand side of (21) equals $-|x|^2(1 - 2\|PBK\|)$ there, and so $1-2\|PBK\|$ cannot be positive since $A$ is not stable. It is also straightforward to show (20) directly: just multiply (10) on the left by $v^T$ and on the right by $v$, where $v$ is a normalized eigenvector of $A^T P + PA$ with a nonnegative eigenvalue. □

### 2.3 Radial and spherical quantization

In the above developments, the required bounds on the quantization error do not depend on the size of the state. This leads to uniform quantization, in the sense that quantization points are distributed uniformly over the region of interest. However, it is well known that more efficient quantization schemes are those which provide lower precision far away from the origin and higher precision close to the origin. Quantizers with a logarithmic scale are particularly useful; see [7]. Loosely speaking, with logarithmic quantization one has the same number of quantization points in the vicinity of every sphere centered at the origin in the state space, whereas with uniform quantization this number grows with the radius. This observation suggests introducing a "direct product" of one quantizer on a unit sphere and another along the radial direction, which is what we do next.

Let us write
$$x = |x|\mathrm{vers}(x)$$
where
$$\mathrm{vers}(x) := \frac{x}{|x|}.$$
We represent the quantizer accordingly as
$$q(x) = q^r(|x|)q^s(\mathrm{vers}(x))$$
where $q^r$ takes $N_1$ positive real values, $q^s$ takes $N_2$ values on or inside the unit sphere, and $N_1$ and $N_2$ are some positive integers such that $N_1 N_2 \leq N$. This means that we introduce two separate quantizers, one for $|x|$ and the other for $\mathrm{vers}(x)$. Quantization points for the latter are allowed to be inside the unit sphere because, as will become clear later, requiring them to lie on the sphere would be too restrictive. (For example, for the case $N_2 = 1$ it is easy to see that one quantization point at the origin achieves a smaller worst-case quantization error over the sphere than any point on the sphere.) The lattice of quantization points for the resulting overall quantizer $q$ is formed by the $N$ pairwise products of values of $q^r$ and $q^s$.

From the triangle inequality and the fact that $|q^s(\mathrm{vers}(x))| \leq 1$ for all $x$ by construction, we obtain
$$\begin{aligned}|q(x) - x| &\leq \left|q^r(|x|)q^s(\mathrm{vers}(x)) - |x|q^s(\mathrm{vers}(x))\right| + \left||x|q^s(\mathrm{vers}(x)) - |x|\mathrm{vers}(x)\right| \\ &\leq \left|q^r(|x|) - |x|\right| + |x|\left|q^s(\mathrm{vers}(x)) - \mathrm{vers}(x)\right| \\ &= |x|\left(\left|\frac{q^r(|x|)}{|x|} - 1\right| + \left|q^s(\mathrm{vers}(x)) - \mathrm{vers}(x)\right|\right)\end{aligned}$$

Using (21), we conclude that the derivative of $V$ along solutions of the closed-loop system satisfies

$$\dot{V} \leq -|x|^2 \left[1 - 2\|PBK\| \left(\left|\frac{q^r(|x|)}{|x|} - 1\right| + |q^s(\text{vers}(x)) - \text{vers}(x)|\right)\right]. \quad (22)$$

Take two numbers $\lambda$ and $\varepsilon$ in the interval $(0,1)$. Then we have $\dot{V} \leq -\varepsilon|x|^2$ whenever the inequalities

$$\left|\frac{q^r(|x|)}{|x|} - 1\right| \leq \frac{\lambda}{2\|PBK\|} \quad (23)$$

and

$$|q^s(\text{vers}(x)) - \text{vers}(x)| \leq \frac{1 - \lambda - \varepsilon}{2\|PBK\|} \quad (24)$$

are satisfied.

Pick a positive number $M$. To handle (23), we take $q^r$ to be a logarithmic quantizer. Define

$$a := 1 - \frac{\lambda}{2\|PBK\|}$$

and

$$b := 1 + \frac{\lambda}{2\|PBK\|}$$

noting from (20) that $a > 0$. Let

$$q^r(s) := (a^i/b^{i-1})M \qquad \text{for } s \in ((a/b)^i M, (a/b)^{i-1} M), \ i \in \{1, \ldots, N_1\} \quad (25)$$

and use arbitrary tie-breaking to define $q^r$ at the endpoints of the above intervals. Then it is easy to check that (23) holds for all $x$ such that

$$(a/b)^{N_1} M \leq |x| \leq M.$$

In view of the condition (24), we introduce the worst-case quantization error on the unit sphere corresponding to $q^s$:

$$\Delta_s := \max_{|x|=1} |q^s(x) - x|. \quad (26)$$

Performing a Lyapunov analysis similar to the one in the previous subsection, we obtain the following result.

**Lemma 3** *Assume that the inequalities*

$$\lambda_{\min}(P) M^2 > \lambda_{\max}(P) (a/b)^{N_1} M^2$$

*and*

$$\Delta_s \leq \frac{1 - \lambda - \varepsilon}{2\|PBK\|} \quad (27)$$

*are satisfied. Then the ellipsoids*

$$\mathcal{R}_1 := \{x \in \mathbb{R}^n : x^T P x \leq \lambda_{\min}(P) M^2\}$$

*and*

$$\mathcal{R}_2 := \{x \in \mathbb{R}^n : x^T P x \leq \lambda_{\max}(P) (a/b)^{N_1} M^2\}$$

*are invariant regions for the system* (11). *Moreover, all solutions of* (11) *that start in the ellipsoid* $\mathcal{R}_1$ *enter the smaller ellipsoid* $\mathcal{R}_2$ *in finite time. An upper bound on this time is*

$$T = \frac{\lambda_{\min}(P) M^2 - \lambda_{\max}(P) (a/b)^{N_1} M^2}{(a/b)^{2N_1} M^2 \varepsilon}.$$

For fixed $N_1$ and $N_2$, the quantity $\Delta_s$ defined by (26) provides a destabilization measure (for $q^s$). When $K$ is given and $\Delta_s$ satisfies the inequality (27) for some $\lambda, \varepsilon \in (0,1)$, we can construct $q^r$ via (25) and compute an ultimate bound on solutions using Lemma 3. Minimizing $\Delta_s$ corresponds to the following optimization problem:

$$\min_{\mathcal{Q}^s, \mathcal{W}^s} \max_{i \in \{1, \ldots, N_2\}} \max_{x \in W_i^s} |q_i^s - x| \tag{28}$$

where $\mathcal{Q}^s = \{q_1^s, \ldots, q_{N_2}^s\}$ is a set of points on or inside the unit sphere and $\mathcal{W}^s = \{W_1^s, \ldots, W_{N_2}^s\}$ is a partition of the unit sphere. An algorithm for solving this problem will be described in Section 3.1. The quantity (28) will not exceed the right-hand side of (27) if $N_2$ is sufficiently large. The values of $N_1$ and $N_2$ which yield the smallest ultimate bound seem to be difficult to compute analytically and in general depend on the stabilizing feedback gain $K$; however, for a given $K$ we only have a finite number of choices for these integers and so can find optimal values by trying all of them. We remark that in the context of the multimedian problem, the idea of spherical coordinates quantization has been exploited before, and in particular the trade-off between the numbers of values for the radial and the spherical directions has been studied; see [24] and the references therein.

It is straightforward to derive similar results using the norm defined by the Lyapunov function, i.e., $\|x\| := \sqrt{x^T P x}$, instead of the Euclidean norm. This gives rise to an optimization problem on an ellipsoid rather than on a sphere.

## 2.4 Radially weighted quantization

The need for logarithmic quantization patterns is evidenced by the fact that it is the ratio $|e|/|x|$, and not the absolute value of the quantization error $|e|$ itself, that needs to be small. This is clear from the formulas (14) and (21). The approach of Section 2.3 leads to an "aligned" logarithmic quantization pattern, in the sense that quantization points on spheres of different radii are aligned along the same radial directions. This is a consequence of the fact that we treated radial and spherical quantization separately, having passed from the formula (21) via the triangle inequality to the formula (22) and then to the two independent inequalities (23) and (24). However, it is not hard to see that non-aligned quantization patterns can achieve better coverage. This suggests proceeding from (21) in a more direct fashion.

To this end, pick two numbers $M > m > 0$ and consider the ellipsoids

$$\mathcal{R}_1 := \{x \in \mathbb{R}^n : x^T P x \leq \lambda_{\min}(P) M^2\}$$

and

$$\mathcal{R}_2 := \{x \in \mathbb{R}^n : x^T P x \leq \lambda_{\max}(P) m^2\}.$$

Define

$$\Delta_{rw} := \max_{x \in \mathcal{R}_1 \setminus \mathcal{R}_2} \frac{|q(x) - x|}{|x|}. \tag{29}$$

The following result then easily follows by virtue of (21).

**Lemma 4** *Assume that*

$$\lambda_{\min}(P) M^2 > \lambda_{\max}(P) m^2 \tag{30}$$

*and*

$$\Delta_{rw} \leq \frac{1 - \varepsilon}{2\|PBK\|} \tag{31}$$

*for some $\varepsilon > 0$. Then the ellipsoids $\mathcal{R}_1$ and $\mathcal{R}_2$ are invariant regions for the system* (11). *Moreover, all solutions of* (11) *that start in the ellipsoid $\mathcal{R}_1$ enter the smaller ellipsoid $\mathcal{R}_2$ in finite time. An upper bound on this time is*

$$T = \frac{\lambda_{\min}(P)M^2 - \lambda_{\max}(P)m^2}{m^2 \varepsilon}.$$

The quantity $\Delta_{rw}$ defined by (29) provides another destabilization measure for $q$, in relation to a pair of numbers $M > m > 0$. Given a stabilizing feedback gain $K$, we can check the inequalities (30) and (31) and, if they are satisfied, obtain an ultimate bound on solutions from Lemma 4. This leads us to the following optimization problem:

$$\min_{\mathcal{Q}, \mathcal{W}} \max_{i \in \{1, \ldots, N\}} \max_{x \in W_i} \frac{|q_i - x|}{|x|} \qquad (32)$$

where $\mathcal{Q} = \{q_1, \ldots, q_N\}$ is a set of quantization points and $\mathcal{W} = \{W_1, \ldots, W_N\}$ is a partition of the annulus $\{x \in \mathbb{R}^n : m < |x| < M\}$ into quantization regions. The optimization problem (32) is different in structure from the ones we encountered earlier, and apparently has not been studied in the locational optimization literature. We henceforth call it the *radially weighted multicenter problem*. It turns out that while this problem is more challenging than the others, it is still computationally tractable. We will develop an algorithm for solving it in Section 3.2. In particular, the inequality (31) will hold for a given $K$ if $N$ is sufficiently large.

## 2.5 Dynamic quantization

So far we have been discussing quantizer design with the understanding that it can only be done once and cannot be changed on-line. The results obtained with this static quantization approach assert the existence of two nested invariant regions such that all trajectories of the quantized system starting in the bigger region approach the smaller one. Let us now suppose that we can recompute the locations of quantization regions and quantization points on-line. Then it is actually possible to achieve global asymptotic stability of the quantized closed-loop system. In this section we briefly explain how this can be done. The discussion closely parallels that of [2, 12], and we refer the reader to these papers for details. In fact, the quantization schemes considered here all fit into the general framework developed in [12].

A control strategy which incorporates dynamic quantization consists of two stages. The objective of the first, "zooming-out" stage is to find a number $M$ such that $x(t_0)$ belongs to the set $\mathcal{R}_1$ defined by (15) for linear systems or (5) for nonlinear systems. This is easy to accomplish, for example, by using a binary quantizer which takes one value inside the ball of a variable radius $\mu$ centered at the origin and another value outside this ball. If we set the control to zero and increase $\mu$ fast enough to dominate the expansiveness of the open-loop system, then at some time $t_0$ we will know that the state lies inside the ball of radius $\mu(t_0)$, at which point we can set $M := \sqrt{\lambda_{\max}(P)/\lambda_{\min}(P)}\mu(t_0)$ in the linear case or $M := \alpha_1^{-1} \circ \alpha_2(\mu(t_0))$ in the nonlinear case.

During the second, "zooming-in" stage, we achieve asymptotic convergence to the origin by applying one of the results expressed by Lemmas 1, 2, 3, and 4, depending on which situation is being considered. We know that if an appropriate inequality involving the size of the quantization error is satisfied, then $x(t_0 + T) \in \mathcal{R}_2$. At time $t_0 + T$, we redefine $M$ so that $\mathcal{R}_2$ coincides with $\mathcal{R}_1$ for the new value of $M$. Recomputing the quantizer for this smaller region and repeating the procedure, we have $x(t) \to 0$ as $t \to \infty$ (for linear systems, the convergence is exponential); one can also show that the origin is Lyapunov stable. This technique effectively mimics a logarithmic quantizer with a countable set of values, of the kind considered in [7]. It is important to recognize that since $\mathcal{R}_2$ is an invariant region for the closed-loop system, it is acceptable to use a conservative upper bound on $T$.

When the quantization is static, our goal is to make the attracting invariant region $\mathcal{R}_2$ as small as possible. With dynamic quantization, we only need to ensure that $\mathcal{R}_2$ is a strict subset of $\mathcal{R}_1$. (In fact, since convergence may become slower closer to the boundary of $\mathcal{R}_2$, we might want to make the annulus $\mathcal{R}_1 \setminus \mathcal{R}_2$ as small as possible, so that we can achieve better coverage on it with a given number of quantization points and obtain faster convergence.) The optimization problems formulated above remain relevant, except that by passing from static to dynamic quantization we basically pass from optimal to suboptimal quantizer design. In the context of Section 2.2, this corresponds to simply making sure that $\Delta$ is small enough so that the inequality (17) is satisfied. In Section 2.3 we just need to enforce the inequality (27) and can choose a small value for $N_1$; for example, we can let $N_1 = 1$, thereby eliminating radial quantization altogether. In Section 2.4 we do not need $m$ to be small compared to $M$.

**Remark 3** In the above discussion we assumed that quantization regions and quantization points can be recomputed infinitely many times and as often as needed. Of course, having to recompute the quantization parameters too frequently may be undesirable. An interesting direction for future research is to incorporate a cost on the number and frequency of these "switching events" into the quantizer design problem. □

## 3 Continuous multicenter problems in facility location

In this section we present a class of optimization problems related to the field of facility location; see the discussion in Section 1 and the survey [5]. The facility location problem we consider will have as special cases the optimization problems studied in Section 2, and in particular the problems (8) from Section 2.1, (28) from Section 2.3, and (32) from Section 2.4.

Let us review some preliminary concepts. Given a compact connected region $\mathcal{D} \subset \mathbb{R}^n$ and a set of $N$ points $\mathcal{Q} = \{q_1, \ldots, q_N\}$ in $\mathbb{R}^n$, the *Voronoi partition* $\mathcal{V} = \{V_1, \ldots, V_N\}$ of $\mathcal{D}$ generated by $\mathcal{Q}$ is defined according to

$$V_i := \{x \in \mathcal{D} : |x - q_i| \leq |x - q_j| \ \forall j \neq i\}. \tag{33}$$

When it is useful to emphasize the dependency on $\mathcal{Q}$, we shall write $\mathcal{V}(\mathcal{Q})$ or $V_i(\mathcal{Q})$. When $\mathcal{D}$ is a polytope in $\mathbb{R}^n$, each $V_i$ is a polytope, otherwise $V_i$ is the intersection between a polytope and $\mathcal{D}$. The faces of the polytope which defines $V_i$ are given by hyperplanes of points in $\mathbb{R}^n$ that are equidistant from $q_i$ and $q_j$, $j \neq i$; among the latter, only "neighboring" points play a role. Note that this (standard) construction remains valid when $\mathcal{D}$ is a lower-dimensional subset of $\mathbb{R}^n$, such as a sphere. We refer to [4, 18] for comprehensive treatments of Voronoi partitions.

Let $\mathcal{Q} = \{q_1, \ldots, q_N\}$ be a collection of points in $\mathbb{R}^n$ and let $\mathcal{W} = \{W_1, \ldots, W_N\}$ be a partition of $\mathcal{D}$. In what follows, we shall concern ourselves with the function

$$\mathcal{H}(\mathcal{Q}, \mathcal{W}) := \max_{i \in \{1, \ldots, N\}} \max_{x \in W_i} \phi(x) f(|x - q_i|) \tag{34}$$

where $\phi : \mathcal{D} \to [0, \infty)$ is continuous non-negative and $f : [0, \infty) \to [0, \infty)$ is continuous, non-decreasing and unbounded. We also assume $\phi$ does not identically vanish on $\mathcal{D}$. We investigate the optimization problem

$$\min_{\mathcal{Q}, \mathcal{W}} \mathcal{H}(\mathcal{Q}, \mathcal{W}) \tag{35}$$

and refer to it as the *weighted multicenter problem*. In general, $\mathcal{H}$ is a nonlinear non-convex function of the locations $\mathcal{Q}$ and of the partition $\mathcal{W}$. Accordingly, its global minima can be obtained only numerically via nonlinear programming algorithms. However, this and related facility location problems [23, 22, 6] have some peculiar structure that helps us characterize optimal solutions and design useful iterative algorithms. Let us start by considering the *weighted 1-center problem* over $\mathcal{D}$, i.e., assume $N = 1$.

**Lemma 5** *The function $\mathcal{H}_1 : \mathbb{R}^n \to [0, \infty)$ defined by*

$$\mathcal{H}_1(q) := \mathcal{H}(\{q\}, \{\mathcal{D}\}) = \max_{x \in \mathcal{D}} \phi(x) f(|x - q|)$$

*is continuous, radially unbounded, and quasiconvex.*[3] *If $f$ is convex and $\phi$ is constant, then $\mathcal{H}_1$ is convex.*

PROOF. The function $\mathcal{H}_1$ is continuous because it is the maximum of a compact family of continuous functions. Further, $\mathcal{H}_1$ is radially unbounded because for every $x^* \in \mathcal{D}$ such that $\phi(x^*) > 0$ we have $\mathcal{H}_1(q) \geq \phi(x^*) f(|x^* - q|)$ and $f$ is unbounded. To show the other statements, we invoke certain properties of convex and quasiconvex functions; see Sections 2.2 and 2.4 in [1]. At fixed $x$, the function $q \mapsto f(|x - q|)$ is quasiconvex because it is the composition of a convex function with a non-decreasing function. Furthermore, if $f$ is non-decreasing and convex, then $q \mapsto f(|x - q|)$ is convex because, at fixed $x$, it is the composition of a convex function with a convex non-decreasing function. If $f$ is convex, then $q \mapsto \max_{x \in \mathcal{D}} f(|x - q|)$ is convex because it is the pointwise supremum over a set of convex functions. For general $\phi$ and $f$, the function $\mathcal{H}_1$ is quasiconvex because it is the weighted pointwise supremum of quasiconvex functions. □

Next, we let $\text{co}(\mathcal{D})$ denote the convex hull of $\mathcal{D}$ and study global minima of $\mathcal{H}_1$ as follows.

**Lemma 6** *The set of global minimum points for $\mathcal{H}_1$ is compact, convex and has a non-empty intersection with $\text{co}(\mathcal{D})$. If $f$ is strictly increasing, then all global minimum points belong to $\text{co}(\mathcal{D})$.*

PROOF. The fact that the set of global minimum points is compact and convex is an immediate consequence of continuity, radial unboundedness, and quasiconvexity. Let us prove the non-empty intersection with $\text{co}(\mathcal{D})$. Suppose that $q^* \notin \text{co}(\mathcal{D})$ is a global minimum point for $\mathcal{H}_1$. Let $p^* \in \text{co}(\mathcal{D})$ be the closest point to $q^*$, i.e., $p^* := \text{argmin}_{x \in \text{co}(\mathcal{D})} |q^* - x|$. Then $|x - p^*| < |x - q^*|$ for all $x \in \mathcal{D}$, so that, for all $x \in \mathcal{D}$, we have $\phi(x) f(|x - p^*|) \leq \phi(x) f(|x - q^*|) \leq \max_{x \in \mathcal{D}} \phi(x) f(|x - q^*|) = \mathcal{H}_1(q^*)$. Therefore, $\mathcal{H}_1(p^*) = \max_{x \in \mathcal{D}} \phi(x) f(|x - p^*|) \leq \mathcal{H}_1(q^*)$ and $p^*$ also belongs to the set of global minimum points. When $f$ is strictly increasing, the previous argument leads to $\mathcal{H}_1(p^*) < \mathcal{H}_1(q^*)$, which contradicts the assumption that $q^*$ is a global minimum. □

The last two lemmas show that the weighted 1-center problem over $\mathcal{D}$ is a quasiconvex optimization problem, i.e., it consists in minimizing the quasiconvex function $\mathcal{H}_1$ over the convex set $\text{co}(\mathcal{D})$. It is known that every quasiconvex optimization problem can be solved by iterative techniques (via a bisection algorithm solving a convex feasibility problem at each step; see Section 3.2 in [1]). We call $q^*(\mathcal{D})$ a *weighted center* of the region $\mathcal{D}$ if it is a (possibly non-unique) global minimum point:

$$q^*(\mathcal{D}) := \underset{q \in \text{co}(\mathcal{D})}{\text{argmin}} \max_{x \in \mathcal{D}} \phi(x) f(|x - q|).$$

Now, it is useful to consider the weighted multicenter problem again, and define $\mathcal{W} \mapsto \mathcal{Q}^*(\mathcal{W})$ as the map that associates to $\mathcal{W}$ a collection of $N$ (possibly non-unique) global minimum points for the corresponding weighted 1-center problems; in other words, $\mathcal{Q}^*(\{W_1, \ldots, W_N\}) := \{q^*(W_1), \ldots, q^*(W_N)\}$.

---
[3]Recall that a *quasiconvex* function is a function defined on a convex domain and with convex sublevel sets.

Note that these weighted centers are well defined in view of the above discussion since each $W_i$ is compact. Finally, define the *Lloyd map* (or the Lloyd algorithm) $\mathcal{L} : (\mathcal{Q}, \mathcal{W}) \mapsto (\mathcal{Q}', \mathcal{W}')$ where $\mathcal{W}' := \mathcal{V}(\mathcal{Q})$ and $\mathcal{Q}' := \mathcal{Q}^*(\mathcal{W}')$.

**Lemma 7** *At fixed $\mathcal{Q}$, the global minimum of $\mathcal{W} \mapsto \mathcal{H}(\mathcal{Q}, \mathcal{W})$ is achieved at $\mathcal{W} = \mathcal{V}(\mathcal{Q})$. At fixed $\mathcal{W}$, a global minimum of $\mathcal{Q} \mapsto \mathcal{H}(\mathcal{Q}, \mathcal{W})$ is achieved at $\mathcal{Q} = \mathcal{Q}^*(\mathcal{W})$. The Lloyd map is a descent algorithm for the cost function $\mathcal{H}$, i.e., an application of the map is guaranteed not to increase $\mathcal{H}$. Given an initial pair $(\mathcal{Q}_0, \mathcal{W}_0)$, the sequence $\{\mathcal{L}^k(\mathcal{Q}_0, \mathcal{W}_0), \ k \geq 0\}$ approaches the largest set invariant under $\mathcal{L}$ such that $\mathcal{H}(\mathcal{L}(\mathcal{Q}, \mathcal{W})) = \mathcal{H}(\mathcal{Q}, \mathcal{W})$.*

Fixed points of the Lloyd map are pairs $(\mathcal{Q}, \mathcal{W})$ such that $\mathcal{W}$ is the Voronoi partition generated by $\mathcal{Q}$ and at the same time the points in $\mathcal{Q}$ are weighted centers for $\mathcal{W}$. It is an open conjecture that global minima of $\mathcal{H}$ are fixed points of the Lloyd map and that the iteration described in the lemma converges to such fixed points. Nevertheless, the algorithm is of interest to us because it is guaranteed to improve a given quantizer design and provides a good indication as to whether or not $N$ is large enough to achieve the control objective.

The classic Lloyd algorithm is tailored to the continuous multimedian problem as it appears, for example, in the problem of fixed-rate minimum-distorsion quantizer design; see [5, 9]. The classic Lloyd algorithm differs from the one considered here only in the fact that the points in $\mathcal{Q}$ are moved to the centroids—as opposed to the weighted centers—of the respective Voronoi regions. (Centroids are solutions of the 1-median problems.)

Next, we consider the specific settings that arise in the quantizer design problems discussed in the previous section. We study the multicenter problem (8), the spherical multicenter problem (28), and the radially weighted multicenter problem (32). To implement the Lloyd algorithm, two tasks must be carried out repeatedly. One consists in computing the Voronoi partition for a given set of points $\mathcal{Q}$, which is accomplished by the standard procedure described earlier. The other amounts to computing a weighted center for each set $W_i$ in a given partition. Thus for each of the specific multicenter problems studied below, we only need to explain how to solve the corresponding 1-center problem. Some additional remarks on the properties of these particular multicenter problems will also be provided.

## 3.1 Multicenter problem

Let us consider the problem (8) arising in Section 2.1. The domain is a ball centered at the origin or, more generally, an ellipsoid, i.e., $\mathcal{D} = \{x \in \mathbb{R}^n : x^T P x \leq 1\}$ for some positive definite symmetric matrix $P$. Note also that the problem (19) arising in Section 2.2 reduces, via a linear change of coordinates, to the multicenter problem considered here in a lower dimension.

In the problem (8), the weighting function $\phi$ and the performance function $f$ are the identity maps. Under these conditions, we refer to the optimization problem (35) simply as the multicenter problem; see [23, 22]. The multicenter problem can be restated as the problem of covering the region $\mathcal{D}$ with (possibly overlapping) balls of smallest radius. If $\mathcal{B}_1 \subset \mathbb{R}^n$ is the unit ball centered at the origin, and if $R\mathcal{B}_1 + q$ denotes the ball of radius $R$ centered at a point $q$, the problem reads:

$$\min R \qquad \text{subject to} \qquad \bigcup_{i \in \{1,\ldots,N\}} (R\mathcal{B}_1 + q_i) \supseteq \mathcal{D}.$$

Let us analyze the 1-center problem. From Lemma 5 we know that this is a convex optimization problem. For each region $V_i$, the optimal solution $q^*(V_i)$ is the unique center of the minimal-radius enclosing sphere for $V_i$. When $V_i \subset \mathbb{R}^2$ is a polygon, this sphere is referred to as the smallest enclosing circle and algorithms are available to compute it; see [4, Chapter 4]. When $V_i \subset \mathbb{R}^n$ is a polytope, the smallest enclosing ellipsoid (in particular, sphere) can be computed via iterative convex optimization algorithms; see [1]. For a Voronoi region $V_i$ near the boundary of $\mathcal{D}$, which is not a polytope, we can under-approximate it by a polytope generated by the vertices of $V_i$ and suitable additional points on the intersection of $V_i$ with the boundary of $\mathcal{D}$, and then compute the center of this polytope. For a sufficiently close under-approximation, this center will also be the center of $V_i$.

When $\mathcal{D}$ is a unit cube in $\mathbb{R}^n$, the optimal value of the problem (8) satisfies the bounds

$$\frac{1}{2\lfloor \sqrt[n]{N} \rfloor} \leq \Delta \leq \frac{\sqrt{n}}{2\lfloor \sqrt[n]{N} \rfloor}.$$

The upper bound is easily obtained by constructing a uniform cubical quantization pattern, while the lower bound is known as Sukharev's lower bound on dispersion [17, 15]. In the present case when $\mathcal{D}$ is a ball, it is straightforward to obtain similar bounds by considering inscribed and superscribed cubes for $\mathcal{D}$. The upper bound can be used to evaluate the convergence of the Lloyd algorithm. When the lower bound on $\Delta$ is not small enough for the inequality (4) or (17) to hold, it indicates that a different destabilization measure and/or a different stabilizing feedback law must be used.

It is also useful to recall some known facts about the multimedian problem. It is conjectured in [8] that for $N$ sufficiently large, the optimal quantizer with respect to the uniform probability density is given by a tessellation (i.e., translation and rotation) of a fixed polytope, except near the boundary of the region of interest. In two dimensions, polygons that can give rise to such tessellations are equilateral triangles, rectangles, and regular hexagons. Among these, the hexagon is optimal, because it has the smallest mean-square quantization error with respect to its centroid per unit volume. This result remains true if we consider the worst-case rather than mean-square quantization error, which is the quantity being minimized in the multicenter problem. The hexagon achieves the smallest error with respect to its center. (For the unit volume regular hexagon this error is approximately 0.62, compared with 0.707 for the square and 0.936 for the equilateral triangle; the unit-volume disk gives the error of 0.564 but disks cannot be used to obtain tessellations.)

The spherical multicenter problem (28) from Section 2.3 corresponds to the setting where $\mathcal{D} = \{x \in \mathbb{R}^n : |x| = 1\}$ is the unit sphere in $\mathbb{R}^n$. Since the spherical multicenter problem is formulated in terms of the Euclidean distance in $\mathbb{R}^n$, Voronoi partitions of the sphere can be constructed as explained earlier for the general case. Voronoi regions will be *spherical polytopes*, i.e., intersections of polytopes with the unit sphere. The center of a spherical polytope $V_i$ is the center of the minimal-radius enclosing sphere for $V_i$. We can consider a polytope in $\mathbb{R}^n$ generated by the vertices of $V_i$ and perhaps some other points in $V_i$. If enough points are taken, then the center of this polytope will also be the center of $V_i$. As we explained earlier, computing the center of a polytope is a computationally tractable task.

## 3.2 Radially weighted multicenter problem

Here, we study the problem (32) formulated in Section 2.4, namely,

$$\min_{\mathcal{Q}, \mathcal{W}} \max_{i \in \{1,\ldots,N\}} \max_{x \in W_i} \frac{|q_i - x|}{|x|}$$

where the domain is the spherical annulus $\mathcal{D} = \{x \in \mathbb{R}^n : m < |x| < M\}$. We consider the corresponding radially weighted 1-center problem over a set $V \subset \mathcal{D}$:

$$\min_{q \in \operatorname{co}(V)} \max_{x \in V} \frac{|q - x|}{|x|}. \tag{36}$$

The problem is well-posed because $V$ is a subset of $\mathcal{D}$ and therefore does not contain the origin. In what follows, we take $V$ to be a polytope; if it is not, we approximate it by a polytope as before. We begin by making the following observation.

**Lemma 8** *The optimal cost in the problem (36) is smaller than 1 if and only if the set $V$ is separated from the origin by a hyperplane.*

PROOF. Suppose first that $V$ is separated from the origin by a hyperplane, so that $0 \notin \mathrm{co}(V)$. Let $\tilde{q}$ be the projection of the origin onto $\mathrm{co}(V)$, i.e., $\tilde{q} := \mathrm{argmin}_{x \in \mathrm{co}(V)} |x|$. By construction, $|x - \tilde{q}| < |x|$ for all $x \in V$, hence $\max_{x \in V} |\tilde{q} - x|/|x| < 1$. This implies that the optimal cost in the problem (36) is less than 1. To prove the converse, suppose on the contrary that $0 \in \mathrm{co}(V)$. This means that the origin lies on the line segment between two points $x_1, x_2 \in V$. For the optimal cost to be less than 1, the optimal point $q^*$ must belong to the open ball $\{q \in \mathbb{R}^n : |q - x_1| < |x_1|\}$ as well as to the open ball $\{q \in \mathbb{R}^n : |q - x_2| < |x_2|\}$. But the intersection between these two sets is empty, which is a contradiction. □

We shall henceforth assume that the set $V$ is separated from the origin by a hyperplane. For $N$ sufficiently large, the initial quantization points can be chosen in such a way that each of the resulting Voronoi regions indeed has this property. Since the Lloyd algorithm decreases the cost, Lemma 8 implies that all Voronoi regions will then have this property at every step of the iteration.

From Lemma 5 and 6 we know that the problem (36) is quasiconvex and can thus be handled by iterative convex optimization algorithms. However, its structure can be utilized to obtain a solution more constructively. Let us first present an equivalent formulation of this optimization problem.

**Lemma 9** *Let $V$ be a polytope separated from the origin by a hyperplane. Consider the problem of finding the sphere with center $c$ and radius $r$ which encloses $V$ and minimizes $r/|c|$. Let $(c^*, r^*)$ be the parameters of the optimal sphere. Then the optimal value for the problem (36) is $\gamma^* := r^*/|c^*|$ and the optimal point is $q^* := (1 - (\gamma^*)^2)c^*$.*

PROOF. Let $f_q(x) := |q - x|/|x|$. In the problem (36), we search for $q$ that minimizes the value of the function $f_q$ on its smallest level set enclosing $V$. For any $\gamma > 0$, the $\gamma$-level set of $f_q$ is described by

$$|q - x|^2 - \gamma^2 |x|^2 = 0.$$

Because $V$ is separated from the origin by a hyperplane, we know from Lemma 8 that the optimal value of $\gamma$ is smaller than 1. Thus from here on we will only be interested in $\gamma < 1$. A square completion argument leads to

$$\frac{1}{1-\gamma^2}\left(|q-x|^2 - \gamma^2|x|^2\right) = \left|x - \frac{q}{1-\gamma^2}\right|^2 - \frac{\gamma^2}{(1-\gamma^2)^2}|q|^2,$$

so that the $\gamma$-level set of $f_q$ is the sphere $|x - c|^2 = r^2$, with center $c := q/(1 - \gamma^2)$ and radius $r := \gamma|c|$. In the new variables $(c, r)$, we must minimize $\gamma = r/|c|$ among all spheres enclosing $V$. □

Note that the point $q^*$ belongs to $\mathrm{co}(V)$ by Lemma 6, while $c^*$ might not. The above result leads us to considering the problem

$$\min_{c \in \mathbb{R}^n, r \in \mathbb{R}} \quad \gamma^2(c, r) := \frac{r^2}{|c|^2} \quad \text{subject to} \quad |c - v_i|^2 \leq r^2, \quad i = 1, \ldots, p \qquad (37)$$

where $v_1, \ldots, v_p$ are the vertices of the polytope $V$. This is an optimization problem subject to inequality constraints, which can be solved with a finite number of computations. The idea is to enumerate active constraints, according to the procedure described in the following algorithm:

1: **for** all subsets $S$ of the set of vertices of $V$ **do**
2:   compute the $(c_S, r_S)$-sphere minimizing $\gamma^2$ among all $(c, r)$-spheres touching all points in $S$
3: **end for**
4: discard all $(c_S, r_S)$-spheres not containing all vertices of $V$
5: find global minima for (37) by comparing the values of $r_S^2/|c_S|^2$ among all remaining candidate spheres

Steps 4 and 5 are straightforward comparison checks. Regarding step 1, it turns out we can restrict our search to sets $S$ containing at least two vertices of $V$, by virtue of the following result.

**Lemma 10** *The optimal sphere for the problem* (37) *touches at least two vertices of $V$, i.e., at least two constraints are active at the minimum.*

PROOF. The proof is by contradiction. Suppose that the optimal sphere touches only one vertex. We denote this vertex by $v$ and assume, performing an affine coordinate change, that it has coordinates $(1, 0, \ldots, 0)^T$. Let $c = (\bar{x}_1, \bar{x}_2, \ldots, \bar{x}_n)^T$. Then we are led to minimizing

$$\gamma^2(\bar{x}_1, \bar{x}_2, \ldots, \bar{x}_n) = \frac{(\bar{x}_1 - 1)^2 + \sum_{i=2}^n \bar{x}_i^2}{\sum_{i=1}^n \bar{x}_i^2} = 1 + \frac{1 - 2\bar{x}_1}{\sum_{i=1}^n \bar{x}_i^2}. \tag{38}$$

Let us show that this function has no critical points besides the pole at the origin and the zero at $v$. We have

$$\frac{\partial \gamma^2}{\partial \bar{x}_1} = \frac{2\left(\bar{x}_1^2 - \sum_{i=2}^n \bar{x}_i^2 - \bar{x}_1\right)}{\left(\sum_{i=1}^n \bar{x}_i^2\right)^2} \tag{39}$$

and

$$\frac{\partial \gamma^2}{\partial \bar{x}_i} = \frac{-2\bar{x}_i(1 - 2\bar{x}_1)}{\left(\sum_{i=1}^n \bar{x}_i^2\right)^2}, \qquad i \neq 1. \tag{40}$$

In view of the formula (40), every critical point satisfies either $\bar{x}_1 = 1/2$ or $\bar{x}_i = 0$ for all $i \neq 1$. In the first case, the formula (39) implies that we must have $\sum_{i=2}^n \bar{x}_i^2 = -1/2$, and this equation has no solution. In the second case, (39) gives two solutions: $\bar{x}_1 = 0$ (pole at 0) and $\bar{x}_1 = 1$ (zero at $v$). The pole at the origin is not a minimum. The zero at $v$ corresponds to the sphere of radius 0 centered at $v$, which is not a feasible solution because it does not enclose $V$. In summary, we have shown that the optimal sphere cannot touch only a single vertex of $V$. □

Regarding step 2, we need to minimize $\gamma^2$ over spheres passing through two or more vertices of $V$. Spheres passing through $l$ generic points in $\mathbb{R}^n$ are parameterized by $n + 1 - l$ variables. A convenient parameterization is obtained by intersecting hyperplanes of points equidistant from pairs of points from a given set. Coordinates of the points on the intersection are given by affine functions of $n + 1 - l$ free parameters. Note that the radius $r$ of the sphere is uniquely determined by its center $c$ and the vertices of $V$ which lie on the sphere. It is straightforward to verify that the function $\gamma^2$ in (37) is a rational function whose numerator and denominator are quadratic inhomogeneous polynomials in these free parameters, and that critical points of $\gamma^2$ are solutions of $n + 1 - l$ quadratic equations in the same number of unknowns. According to Bezout's theorem, this genericaly gives $2^{n+1-l}$ candidate optimal spheres (see [3]). Step 2 is completed by choosing the one with the smallest radius.

As an example of step 2, let us work out the planar case. When $n = 2$, the problem reduces to finding critical points of $\gamma$ for circles passing through $l$ vertices of $V$, where $l > 1$ by Lemma 8. Since for $l > 2$ there is at most one circle passing through the corresponding vertices, we only need to explain how to solve this problem for $l = 2$. For convenience, let us consider an affine change of coordinates which places the two vertices at $(1, 0)^T$ and $(-1, 0)^T$ and the origin at some point $(x_0, y_0)^T$. Without loss of generality, assume that $y_0 \geq 0$. The center of the circle is denoted by $c = (\bar{x}, \bar{y})^T$. We know that $c$ must be equidistant from the two vertices, hence $\bar{x} = 0$. Then we have

$$\gamma^2 = \frac{1 + \bar{y}^2}{x_0^2 + (\bar{y} - y_0)^2}$$

and so

$$\frac{\partial \gamma^2}{\partial \bar{y}} = \frac{2\bar{y}(x_0^2 + (\bar{y} - y_0)^2) - 2(\bar{y} - y_0)(1 + \bar{y}^2)}{(x_0^2 + (\bar{y} - y_0)^2)^2} = \frac{-2\bar{y}^2 y_0 + 2\bar{y}(x_0^2 + y_0^2 - 1) + 2y_0}{(x_0^2 + (\bar{y} - y_0)^2)^2}.$$

Equating the numerator to 0, we arrive at the equation

$$-\bar{y}^2 y_0 + \bar{y}(x_0^2 + y_0^2 - 1) + y_0 = 0.$$

In the special case when $y_0 = 0$, this reduces to $\bar{y}(x_0^2 - 1) = 0$. Since $x_0 = \pm 1$ corresponds to one of the vertices being at the origin, which cannot happen by our earlier assumption, the solution is $\bar{y} = 0$ (as is also clear from symmetry). When $y_0 \neq 0$, the minimum is achieved at

$$\bar{y} = \frac{x_0^2 + y_0^2 - 1 - \sqrt{(x_0^2 + y_0^2 - 1)^2 + 4y_0^2}}{2y_0} < 0.$$

(Note that this goes to 0 as $y_0$ approaches 0 or $\infty$.)

# References

[1] S. Boyd and L. Vandenberghe. Convex optimization. Preprint, December 2001.

[2] R. W. Brockett and D. Liberzon. Quantized feedback stabilization of linear systems. *IEEE Trans. Automat. Control*, 45:1279–1289, 2000.

[3] J. L. Coolidge. *A Treatise on Algebraic Plane Curves*. Dover, New York, 1959.

[4] M. de Berg, M. van Kreveld, and M. Overmars. *Computational Geometry: Algorithms and Applications*. Springer Verlag, New York, NY, 1997.

[5] Z. Drezner, editor. *Facility Location: A Survey of Applications and Methods*. Springer Series in Operations Research. Springer Verlag, New York, NY, 1995.

[6] Q. Du, V. Faber, and M. Gunzburger. Centroidal Voronoi tessellations: applications and algorithms. *SIAM Review*, 41(4):637–676, 1999.

[7] N. Elia and S. K. Mitter. Stabilization of linear systems with limited information. *IEEE Trans. Automat. Control*, 46:1384–1400, 2001.

[8] A. Gersho. Asymptotically optimal block quantization. *IEEE Trans. Information Theory*, 25:373–380, 1979.

[9] R. M. Gray and D. L. Neuhoff. Quantization. *IEEE Transactions on Information Theory*, 44(6):2325–2383, 1998. Commemorative Issue 1948-1998.

[10] H. Ishii and B. A. Francis. Stabilizing a linear system by switching control with dwell time. *IEEE Trans. Automat. Control*, 2002. To appear.

[11] D. Liberzon. Nonlinear stabilization by hybrid quantized feedback. In N. Lynch and B. H. Krogh, editors, *Proc. Third International Workshop on Hybrid Systems: Computation and Control*, volume 1790 of *Lecture Notes in Computer Science*, pages 243–257. Springer, 2000.

[12] D. Liberzon. Hybrid feedback stabilization of systems with quantized signals. *Automatica*, 2001. Submitted.

[13] S. P. Lloyd. Least squares quantization in PCM. *IEEE Transactions on Information Theory*, 28(2):129–137, 1982. Presented as Bell Laboratory Technical Memorandum at a 1957 Institute for Mathematical Statistics meeting.

[14] J. Lunze, B. Nixdorf, and J. Schröder. Deterministic discrete-event representations of linear continuous-variable systems. *Automatica*, 35:395–406, 1999.

[15] J. Matousek. *Geometric Discrepancy: An Illustrated Guide.* Number 18 in Algorithms and Combinatorics. Springer Verlag, New York, NY, 1999.

[16] G. N. Nair and R. J. Evans. Exponential stabilisability of multidimensional linear systems with limited data rates. *Automatica*. To appear.

[17] H. Niederreiter. *Random Number Generation and Quasi-Monte Carlo Methods.* Number 63 in Cbms-Nsf Regional Conference Series in Applied Mathematics. Society for Industrial & Applied Mathematics, 1992.

[18] A. Okabe, B. Boots, K. Sugihara, and S. N. Chiu. *Spatial Tessellations: Concepts and Applications of Voronoi Diagrams.* Wiley Series in Probability and Statistics. John Wiley & Sons, New York, NY, second edition, 2000.

[19] A. Okabe and A. Suzuki. Locational optimization problems solved through Voronoi diagrams. *European Journal of Operational Research*, 98(3):445–56, 1997.

[20] E. D. Sontag. Smooth stabilization implies coprime factorization. *IEEE Trans. Automat. Control*, 34:435–443, 1989.

[21] E. D. Sontag and Y. Wang. On characterizations of the input-to-state stability property. *Systems Control Lett.*, 24:351–359, 1995.

[22] A. Suzuki and Z. Drezner. The p-center location problem in an area. *Location Science*, 4(1/2):69–82, 1996.

[23] A. Suzuki and A. Okabe. Using Voronoi diagrams. In Drezner [5], pages 103–118.

[24] P. F. Swaszek and J. B. Thomas. Multidimensional spherical coordinates quantization. *IEEE Transactions on Information Theory*, 29(4):570–6, 1983.